\documentclass [leqno,12pt]{article}
\usepackage[utf8]{inputenc} 
\usepackage[T1]{fontenc} 
\usepackage{amsmath}
\usepackage{amsfonts}
\usepackage{amsthm}
\usepackage{color}
\usepackage{graphicx}
\usepackage{epstopdf}

\usepackage[small,nohug,heads=vee]{diagrams}
\diagramstyle[labelstyle=\scriptstyle]

\newlength{\oldparindent}
\newcommand{\un}{{\bf 1}}

\newcommand{\cL}{{\mathbb {L}}}

\newcommand{\bpf}{\begin{preuve}}
\newcommand{\epf}{ \end{preuve} \medskip}

\newcommand{\benum}{\begin{enumerate}}
\newcommand{\eenum}{\end{enumerate}}

\newcommand{\bitem}{\begin{itemize}}
\newcommand{\eitem}{\end{itemize}}

\newcommand{\brmq}{\begin{rmq}}
\newcommand{\ermq}{\end{rmq}}

\newcommand{\brmqs}{\begin{rmqs}}
\newcommand{\ermqs}{\end{rmqs}}

\newcommand{\bapp}{\begin{application}}
\newcommand{\eapp}{\end{application}}

\newcommand{\bapps}{\begin{applications}}
\newcommand{\eapps}{\end{applications}}

\newcommand{\bdefi}{\begin{definition}}
\newcommand{\edefi}{\end{definition}}

\newcommand{\beq}{\begin{equation}}
\newcommand{\eeq}{\end{equation}}

\def\bpm{\begin{pmatrix}}
\def\epm{\end{pmatrix}}

\newcommand{\bcas}{\begin{cases}}
\newcommand{\ecas}{\end{cases}}

\newcommand{\bex}{\begin{exemp}}
\newcommand{\eex}{\end{exemp}}

\newcommand{\bexs}{\begin{exemps}}
\newcommand{\eexs}{\end{exemps}}

\newcommand{\bprop}{\begin{proposition}}
\newcommand{\eprop}{\end{proposition}}

\newcommand{\bthm}{\begin{theoreme}}
\newcommand{\ethm}{\end{theoreme}}

\newcommand{\bcor}{\begin{corollaire}}
\newcommand{\ecor}{\end{corollaire}}

\newcommand{\blem}{\begin{lemme}}
\newcommand{\elem}{\end{lemme}}

\newcommand{\beqna}{\begin{eqnarray}}
\newcommand{\eeqna}{\end{eqnarray}}

\newcommand{\beqnas}{\begin{eqnarray*}}
\newcommand{\eeqnas}{\end{eqnarray*}}

\newcommand{\cA}{{\mathcal A}}

\definecolor{green}{rgb}{0,.7,.2}
\definecolor{orange}{rgb}{0.9,.5,0}

\newcommand{\LL}{{\rm L}}

\def\I{1\!{\rm l}}
\def\tr{\textmd{trace}\,}

\def\lag{\langle}
\def\rag{\rangle}

\def\det{{ \rm{det}}}  

\def\dis{\displaystyle }

\def\Id{{\rm{Id}}} 

\def\eps{\varepsilon}

\def\cA{{\mathcal A }}
\def\cB{{\mathcal B }}

\def\cD{{\mathcal D}}

\def\cL{{\mathcal L }}

\def\cP{{\mathcal P }}

\def\cS{{\mathcal  S}}

\def\cV{{\mathcal V}}

\def\cX{{\mathcal X}}

\def\bbC{{\mathbb{C}}}

\def\bbE{{\mathbb{E}}}

\newcommand{\bbN}{{\mathbb {N}}}
\newcommand{\bbR}{{\mathbb {R}}}

\newtheorem{theoreme}{Theorem}[section]
\newtheorem{lemme}[theoreme]{Lemma}
\newtheorem{definition}[theoreme]{Definition}
\newtheorem{proposition}[theoreme]{Proposition}
\newtheorem{corollaire}[theoreme]{Corollary}

\newenvironment{exemp}{\noindent{\bf Example. --- }}{\par}
\newenvironment{exemps}{\noindent{\bf Examples}\benum}{\eenum\par}
\newtheorem{rmq}[theoreme]{Remark}
\newtheorem{rmqs}[theoreme]{Remarks}

\newenvironment{preuve}{\noindent{\it Proof. --- }}
{\hfill\rule{1.3mm}{2mm}\par} 
\newenvironment{application}{\noindent{\bf Application. --- }}{\par}
\newenvironment{applications}{\noindent{\bf Applications. --- 
}\benum}{\eenum\par}

\theoremstyle{definition}

\author{Dominique Bakry, Olfa Zribi\\ {\small Institut de Mathématiques de Toulouse,} \\ {\small Universit\'e Paul Sabatier,} \\ {\small
118 
route de Narbonne,} \\ {\small 31062 Toulouse,} \\ {\small FRANCE}}
\date{}

\date{}

\makeatletter \renewcommand{\@oddfoot}{\sl \small
 \hfil \thepage\hfil }
\renewcommand{\@oddhead}{\sl \small
 \hfil }
 \makeatother

\begin{document}

\title{Hypergroup properties for the deltoid model
}

\maketitle

 \abstract{ We investigate the hypergroup property for the family of orthogonal polynomials associated with the deltoid curve in the plane, related to the $A_2$ root system. This provides also the same property for another family of polynomials related to the $G_2$ root system.}

\section{Introduction}

The hypergroup property is a property shared by some orthonormal bases in probability spaces which allows for the complete description of all Markov sequences and for multiplication formulas 
(see Section~\ref{sec.hyperg}).  It has been established for some  families of orthogonal polynomials in dimension 1, and relates in general to some special structure of the underlying space.  This is a quite powerful tool in many areas, ranging from pure analysis and Lie group  to statistics and computer algorithms.  Gasper's theorem establishes this property for  dissymetric Jacobi polynomials   in dimension 1 (see \cite{Gasper70, Gasper71, Gasper72, Koorn75}), and  is an extension of an earlier  result due to Bochner concerning the symmetric case~\cite{Bochner55, Bochner79}. Although many authors revisited this result, a very elegant proof of Gasper's theorem was recently proposed by Carlen, Geronimo and Loss  in~\cite{CGS2010}, arising from the study of Kac's model in statistical mechanics. This proof relies on the construction of an auxiliary model, and may be extended in many other situations. The aim of this paper is to show how this technique applies in the particular case of one of the 2 dimensional extension of Jacobi polynomials, namely  orthogonal polynomials on the deltoid  domain (defined below in Section~\ref{sec.delt.mod}) and associated to the $A_2$ root system.

Jacobi polynomials is the unique family of orthogonal polynomials in dimension 1 (together with their scaled limits Hermite and Laguerre polynomials) which are at the same time  eigenvectors of a second order differential operator, and more precisely of a diffusion operator~\cite{Mazet97}. Although there is no reason a priori for these two properties to be related, it is worth to study the hypergroup property for other families of  such orthogonal polynomials in higher dimension. Among these orthogonal polynomials associated with diffusion operators, some of them arise from root systems,  the Heckman-Opdam Jacobi polynomials~\cite{HeckmanOpdam87, HeckmanOpdam97}. It is an open and challenging problem to analyse this question  for these families in particular.     Partial results have been obtained in this direction in~\cite{ RemlingRosler2014,Rosler2010} in the $BC_n$ case. 

We investigate in this paper the family of orthogonal polynomials in dimension 2 related to the deltoid curve and the  $A_2$  and $G_2$ root systems, introduced by \cite{koorn1,koorn2,koorn3},  as one as the special families of orthogonal polynomials in dimension 2 which are at the same time eigenvectors of symmetric diffusion operators (see~\cite{BOZ2013} and  Section~\ref{sec.delt.mod}). The scheme of Carlen-Geronimo-Loss relies in a crucial way on the fact that the polynomials are eigenvectors of some operator, and that the corresponding eigenspaces have dimension one. In the context that we investigate, this last property fails to be true, and this introduces some extra complexity in the formulation of the main hypergroup property result in Section~\ref{sec.hgp.deltoid1}. This difficulty arises from a symmetry invariance in the deltoid model.  To get rid of this difficulty, we may look at simpler forms of the deltoid model, that is consider only functions which are invariant under this symmetry. This leads to investigate a new model, related to the $G_2$ root system, for which the hypergroup property takes the usual form.  

The paper is organized  follows.   In Section~\ref{sec.hyperg}, we present the hypergroup property  and the  elegant approach initiated by Carlen-Geronimo-Loss~\cite{CGS2010} to obtain it.  In Section~\ref{sec.sym.diff}, we introduce the language of symmetric diffusion operators, which be used through the rest of the paper. The deltoid model is described in Section~\ref{sec.delt.mod}, where we give some details about the structure of the eigenspaces, and show the relation with the $A_2$ root system, through some geometric  interpretations  of the associated operators for two distinct values of the parameter.  The Carlen-Geronimo-Loss methods relies on the construction of some other model space (in general in higher dimension), which projects on the model under study and have some extra symmetry. This model (in our case 6-dimensional) is presented in Section~\ref{6.dim.model}. Then we give in Section~\ref{sec.hgp.deltoid1} the hypergoup property for the deltoid model itself,    and finally in Section~\ref{sec.G2}, we introduce the projected model related to the $G_2$ root system   on which the hypergroup property presents a simpler form.   

\section{\label{sec.hyperg} Hypergroup property and Markov sequences}

The world hypergroup had been first introduced in 1959  by H.S Wall~\cite{Wall59}, to generalize the notion of a group, when the product of two elements is the sum of a finite numbers of elements. Different generalizations had been put forward by C. Dunkl,\cite{Dunkl} Jewett~\cite{Jewett}, see also  the exposition book~\cite{BloomHeyer95}. In our context, we shall mostly be concerned by the hypergroup property, as defined in \cite{BH2008}, and which concerns some properties of orthonormal bases in $\cL^2(\mu)$ spaces, where $\mu$ is a probability measure,  and that  we describe below.

Let $(X, \cA, \mu)$ a probability space, for which some orthonormal basis of $\cL^2(\mu)$ is given, which we suppose countable, and therefore ranked as  $(f_0, f_1, \cdots, f_n, \cdots)$. We suppose moreover that $f_0=1$. A  central question which arises quite often is to determine all sequences $(\lambda_n)$ such that, if one defines a linear  operator $K : \cL^2(\mu)\mapsto \cL^2(\mu)$ through $K(f_n)= \lambda_n f_n$, then this operator is Markov.  This means that $K(\un)=\un $ and $K(f)\geq 0$ whenever $f\geq 0$. Of course, the first condition reduces immediately to $\lambda_0=1$ and the difficulty is to check the positivity preserving property.  We call such sequences Markov sequences, and the problem is known as the Markov Sequence Problem (MSP in short).

When $\sum_n \lambda_n^2 <\infty$, the operator $K$ may be represented through the $\cL^2(\mu\otimes\mu)$  kernel $k(x,y)= \sum_n \lambda_n f_n(x)f_n(y)$ and, since all the functions $f_n$ (except $f_0$) satisfy $\int f_n d\mu=0$, the previous series is oscillating and positivity may not be checked directly from this representation. 

The real parameters $\lambda_n$ are the eigenvalues of the symmetric operator $K$, and it is well known (and quite immediate) that, if $K$ is Markov,  they must satisfy $|\lambda_n|\leq 1$. It is  obvious that if $\lambda= (\lambda_n)$ and $\mu= (\mu_n)$ are Markov sequences , the $(\lambda_n\mu_n)$ is a Markov sequence. Moreover, for any $\theta\in [0,1]$, $\theta\lambda+ (1-\theta)\mu= (\theta\lambda_n + (1-\theta)\mu_n)$ is again a Markov sequence, and  the simple limit of Markov sequences is again a Markov sequence. Therefore the question boils down to the determination the extremal points in the set of Markov sequences.  The hypergroup property allows for such description.

The hypergroup property holds when there exists some point $x_0\in X$ such that, for any $x\in X$, the sequence 
$\frac{f_n(x)}{f_n(x_0)}$ is a Markov sequence. Of course, for this to make sense, one needs for example   some topology on $X$ and require the functions $(f_n)$ to be continuous for this topology.

When the hypergroup property holds, then these sequences $\frac{f_n(x)}{f_n(x_0)}$ are automatically the extremal points for the Markov sequence problem (see~\cite{BH2008}) and then, for any Markov sequence $(\lambda_n)$, there exists a probability measure $\nu$ on $\Omega$ such that 
\beq\label{repr.mkv.sqce}\lambda_n = \int_X \frac{f_n(x)}{f_n(x_0)} \nu(dx).\eeq

To understand this representation,  one should  first  extend the operator $K$ by duality as an operator acting on probability measures. Since formally, $\delta_{x_0}$ has a density which may be written as 
$\sum_n f_n(x)f_n(x_0)$, then $K(\delta_{x_0})$ may be written as $\big(\sum_n \lambda_n f_n(x)f_n(x_0)\big) d\mu(x)$ and, with $\nu(dx)= K(\delta_{x_0})$, one has 
$$\int \frac{f_n(x)}{f_n(x_0)} \nu(dx)= \lambda_n.$$

This of course is not really meaningful beyond the case of finite sets, and one should replace $\delta_{x_0}$ by a smooth approximation of it. It what follows, we shall have a symmetric diffusion operator $ \LL$ with eigenvectors $f_n$, and the associated heat kernel   $P_t= e^{t\LL}$ will be such that $P_t(\delta_{x})$ has a bounded density with respect to $\mu$ for any $t>0$ and any $x$. One may then replace $\delta_{x_0}$ by $P_t(\delta_{x_0})$ in the previous argument, and let then $t\to 0$ to get the representation

When the space $X$ is  a finite set, is is a consequence of the definition that this point $x_0$ must be the point with minimal measure (see~\cite{BH2008}). Moreover, provided we chose the signs of $(f_n)$ such that $f_n(x_0)\geq 0$, it is  immediate (and does not require the finiteness of the space $X$) that the functions $f_n$ reach their maximum (and the maximum value of their modulus) at $x_0$.

This hypergroup property is strongly related with product formulae (see~\cite{KoornSchw}). Indeed, assume that for any $x\in \Omega$, the Markov kernel  $K_x$ with eigenvalues $\frac{f_n(x)}{f_n(x_0)}$ has a $\cL^2$ integrable density
$$K_x(f)(y) = \int f(z) k(x,y,z) \mu(dz).$$
Then, it is readily seen that 
$$K(x,y,z) = \sum_n \frac{f_n(x)f_n(y)f_n(z)}{f_n(x_0)}$$ is a symmetric function of $(x,y,z)$ and we get a product formula
\beq\label{eq.product.formula}f_n(x)f_n(y) =f_n(x_0)  \int f_n(z) K(x,y,z) \mu(dz).\eeq

In \cite{CGS2010}, the authors provide a new and elegant proof of Gasper's result, through a method which proves to be efficient in many other situations, and that we shall describe now. For this, we first require some topology on $X$ and the functions $(f_n)$ to be continuous for this topology.

We assume that there exists a self-adjoint operator  (in general unbounded and defined on a dense subset $\cD\subset \cL^2(\mu)$) $\LL : \cD \mapsto \cL^2(\mu))$ for which the functions $f_n$ are eigenvectors,
that is $ \LL f_n= \mu_n f_n$. Moreover, we assume that the eigenspaces of $\LL$ for the eigenvalues $\mu_n$ are simple.  In our context, one may chose $\cD$ as the set of the finite linear combinations of the functions $f_i$.

 We also assume the existence of an auxiliary  probability space $(Y, \cB, \nu)$ with a self adjoint operator $\hat \LL$, again defined on a  subset $\hat \cD \subset  \cL^2(\nu)$,  $\hat \LL :  \hat \cD \mapsto \cL^2(\nu)$, together with some application $\pi : Y\mapsto X$. For a function $f : X\mapsto \bbR$, we denote by $\pi(f) : Y\mapsto \bbR$, the function $\pi(f)(y)= f(\pi(y))$. We suppose  that $\pi(\cD)\subset \hat \cD$ and that $\hat \LL\pi= \pi \LL$.  This property is  often described as   " $\LL$ is the image of $\hat \LL$ under $\pi$" , or that "$\hat \LL$ projects onto $\LL$ under $\pi$", or  even that "$\pi$ intertwines $\LL$ and $\hat \LL$".  
 
 Moreover, one requires some measurable transformation $\Phi : Y\mapsto Y$ which commutes with the action of $\hat \LL$. That is, once again denoting $\Phi(f)(y)= f(\Phi(y))$, we assume that $\Phi : \hat \cD\mapsto \hat \cD$ and $\hat \LL \Phi= \Phi\hat \LL$.
 
 In order for the next proposition to make sense, we shall impose some topology on  $Y$, and require  $ \pi$ and  $\Phi$ to be continuous.
 
Let $\xi$ be a random variable with values in $Y$, distributed according to $\nu$. From the hypotheses, if is quite clear that the laws of $\pi(\xi)$ and of $\pi(\Phi(\xi))$ are $\mu$.  We look at the conditional law  $k(x,dy)$ of $\pi(\Phi(\xi))$ given $\pi(\xi)$, where $x\mapsto k(x,dy)$ is continuous for the weak topology.
 
 \bprop\label{prop.hpgp.gal}  Assume that there exist some point $x_0\in X$ such that the conditional law $k(x_0, dy)$ is the Dirac mass at some point $x_1\in X$. Then,
 the sequence $\lambda_n=\frac{f_n(x_1)}{f_n(x_0)}$ is a Markov sequence.

 \eprop
 
 \bpf 
 
 We consider the correlation operator, defined as 
 $$K(f)(x)= \int f(y) k(x,dy).$$ It is a Markov operator by construction.  The assumption tells us that, for $f$ continuous, $K(f)(x_0)= f(x_1)$. We shall see that  this operator $K$  corresponds to the choice of the Markov sequence $(\mu_n)= \frac{f_n(x_1)}{f_n(x_0)}$.
 By definition of the conditional expectation, for any pair of $\cL^2(\mu)$ functions
 $$\int_X  f K(g) d\mu= \int_Y f(\pi(y)) g(\pi(\Phi(y)) \nu(dy)= \int_Y \pi (f) \Phi\pi (g) d\nu.$$
 We want to show first that for any $n$, $f_n$ is an eigenvector for $K$. For this, chose any $p$ and consider
 \beqnas \int_X \LL f_p K(f_n))&=&  \int_Y \pi  \LL (f_p) \Phi\pi(f_n)d\nu= \int_Y \hat \LL \pi(f_p) \Phi\pi(f_n)) d\nu\\&=&
 \int_Y \pi (f_p) \hat\LL\Phi\pi (f_n)d\nu= \int_Y \pi (f_p) \Phi\pi \LL(f_n)d\nu\\&=& \lambda_n\int_Y\pi( f_p)  \Phi\pi(f_n) d\nu= \lambda_n \int_X f_p K(f_n) d\mu.
 \eeqnas
 Therefore, this being valid for any $p$, $K(f_n)$ is an eigenvector for $\LL$ with eigenvalue $\lambda_n$. Since the eigenspaces are one dimensional, we deduce that 
 $K(f_n) = \mu_n f_n$. Applying this at the point $x_0$, we see  that $\mu_n = \frac{f_n(x_1)}{f_n(x_0)}$.

 \epf

 If we have at disposal a full family of such transformations $\Phi$ such that the associated points $x_1$(depending on $\Phi$) cover $X$, we conclude to the semigroup property.  The main challenge for proving the hypergroup property,  when we have a basis given as eigenvectors for some operator $\LL$,  is then to construct the space $Y$, the operator $\hat \LL$, the projection $\pi$ and the transformations $\Phi : Y\mapsto Y$, which satisfy the required properties. This is what we are going to do for the deltoid model in Section~\ref{6.dim.model}.
 
  \section{Symmetric diffusion operators\label{sec.sym.diff}}
 We briefly recall in this Section the context of symmetric diffusion operators, following~\cite{bglbook}, a specific context adapted to our setting.
 
For a given probability  space $(X, \cX, \mu)$, we suppose given an algebra $\cA$ of functions  such that 
$\cA\subset \cap_{1\leq p<\infty} \cL^p(\mu)$, where $\cA$ is dense in $\cL^2(\mu)$. A bilinear application   
$\Gamma :~\cA\times \cA\mapsto \cA$ is given such that, $\forall f\in \cA$,  $\Gamma(f,f)\geq 0$.   If $\Phi : \bbR^n\mapsto \bbR$ is  a smooth function such that for any $(f_1, \cdots, f_n)\in \cA^n$, $\Phi(f_1, \cdots, f_n)\in \cA$, then 
$$\Gamma(\Phi(f_1, \cdots, f_k), g)= \sum_i \partial_i \Phi \Gamma(f_i,g).$$ A linear operator $\LL$ is defined through
\beq\label{ipp}\int_X f\LL(g) d\mu= -\int_X \Gamma(f,g) d\mu\eeq and we assume that $\LL$ maps $\cA$ into $\cA$. We extend $\LL$ into a self adjoint operator and we suppose and that $\cA$ is dense in the domain of $\LL$.

Then, for $f= (f_1, \cdots, f_k)$, and again if $\Phi$ is smooth with $\Phi(f_1, \cdots, f_n)$ in $\cA$  whenever $f_i\in \cA$, then 
\beq\label{diff.prop}\LL(\Phi(f))= \sum_1^k \partial_i \Phi(f) \LL(f_i) + \sum_{i,j=1}^k \partial^2_{ij} \Phi(f) \Gamma(f_i,f_j).\eeq
We have from~\eqref{diff.prop} 
$$\Gamma(f,g)= \frac{1}{2} \Big(\LL(fg)-f\LL(g)-g\LL(f)\Big).$$

If $X$ is an open domain in $\bbR^d$, of some smooth manifold, then, in a local system of coordinates $(x^1, \cdots, x^d)$, $\LL$ may in general  be written as
\beq\label{L.loc.coord}\LL(f)= \sum_{ij} g^{ij}(x)\partial^2_{ij}f + \sum_i b^i(x)\partial_i f,\eeq 
while 
\beq\label{Gamma.Loc.coord}\Gamma(f,g)= \sum_{ij} g^{ij}(x) \partial_i f\partial_j g,\eeq
where $g^{ij}= \Gamma(x^i, x^j)$ and $b^i= \LL(x^i)$.
The non negativity of $\Gamma$ translates into the fact that the symmetric matrix $(g^{ij})(x)$ is  everywhere non negative.
Moreover,  $g^{ij}(x)= \Gamma(x^i,x^j)$ and $b^i(x)= \LL(x^i)$. 
Observe then that in order to describe $ \LL$, we just have to describe $ \LL(x^i)$ and $\Gamma(x^i,x^j)$.

When the measure $\mu$ has a positive density $\rho$ with respect to the Lebesgue measure,  then the symmetry property of $ \LL$ translates into
\beq\label{L.sym.rho} \LL= \frac{1}{\rho} \sum_{ij} \partial_i\Big( g^{ij} \rho \partial_j\Big),\eeq
which shows that, in formula~\eqref{L.loc.coord},
\beq\label{eq.bi} b_i = \sum_j \partial_jg^{ij} + g^{ij}\partial_j \log \rho,\eeq
and this formula may be applied in many circumstance to identify the  measure density $\rho$ up to some normalizing constant. Since we shall mainly use this setting for finite measure, we may always assume with no loss of generality that $\mu$ is a probability measure.

A central question is to determine on which set of functions one apply the operator $ \LL$, particularly when $\Omega$ is bounded.  Indeed, this requires to look at some self adjoint extension of $ \LL$, and one needs in general to describe an algebra of functions which is dense in the domain of $ \LL$ and on which  the integration by parts formula~\eqref{ipp} holds true. This is done in general by prescribing boundary conditions on the functions $f$, such as Neumann or Dirichlet conditions.

In our setting however, we will always work on bounded open sets $\Omega\subset \bbR^d$ with piecewise smooth boundary $\partial\Omega$. Our functions $g^{ij}$ and $b^i$ will be smooth in $\Omega$. Moreover, $ g^{ij}(x)$ will be defined and smooth in some neighborhood of $\Omega$.    Then, suppose that in a neighborhood $\cV$ of any regular point $x$ of the boundary, the boundary may be described through $\{F=0\}$, where $F$ is a smooth function, defined in $\cV$ and with real values. Then, our fundamental assumption is that 
\beq\label{eq.bord} \Gamma(F, x_i)=0 \hbox{ on $\partial \Omega\cap \cV$}.\eeq

When this happens, we may chose for $\cA$ the algebra of smooth compactly supported functions defined in a neighborhood of $\Omega$, referred in what follows as "smooth functions", and the integration by parts formula~\eqref{ipp}  holds true  for those functions. In other words, for such operators, there is no need to consider boundary conditions of the functions in~$\cA$ (see~\cite{BOZ2013}). In the context of orthogonal polynomials on bounded domains that we shall consider, this property is always satisfied (see~\cite{BOZ2013}).

We shall make a strong use of  the notion of image operator, to fit with the setting described in Section~\ref{sec.hyperg}.  Suppose that we have a set of functions 
$X= (X_1, \cdots, X_k)$ for which 
$$ \LL(X^i)= B^i(X), \Gamma(X^i, X^j)= G^{ij}(X),$$ then, of any  smooth function $\Phi : \bbR^k\mapsto \bbR$, and  thanks to equation~\eqref{diff.prop}, one has
$$ \LL(\Phi(X))=  \LL_1(\Phi)(X),$$ where 
$$ \LL_1= \sum_{ij} G^{ij} \partial^2_{ij} + \sum_iB^i\partial_i.$$
This is again a symmetric  diffusion operator, defined on the image $\Omega_1= X(\Omega)$, and its reversible measure is the image of $\mu$ through $X$. This new diffusion operator $\LL_1$ is therefore the image of $\LL$ under $\Phi$, in the sense described in Section~\ref{sec.hyperg} : $\Phi \LL_1= \LL \Phi$ by construction.

In the next sections, we shall always work on polynomials, moreover in even dimensions $2k$.  We shall suppose that the coordinates are  paired as $(x_p,y_p), p= 1, \cdots, k$. In this context, it is often quite simpler to use complex coordinates (that is identify $\bbR^{2k}\simeq \bbC^k$), setting $z_p= x_p+iy_p, \bar z_p= x_p-iy_p$ and one has to describe then, using linearity and bilinearity 
$\Gamma(z_p,z_q), \Gamma(\bar z_p, \bar z_q), \Gamma(z_p, \bar z_q)$, together with $ \LL(z_p)$ and $ \LL(\bar z_p)$. 

For example, $\Gamma(z_p,z_p)= \Gamma(x_p,x_p)-\Gamma(y_p,y_p)+ 2i \Gamma(x_p, y_p)$.

The positivity of the metric here may be checked  according to the parity of $k$. For example, a careful inspection shows that indeed, the determinant of the metric in the variables $(z_p,\bar z_p)$ is $(-1)^k 4^k\det(g)$, where the  determinant   computed in real variables. 
 
 \section{The deltoid model\label{sec.delt.mod}}
 
 The deltoid curve is a degree 4 algebraic plane curve which may be parametrized as 
\begin{equation*}
x(t) =\frac{1}{3}( 2\cos  t + \cos 2t), \quad y(t) = \frac{1}{3}(2\sin t - \sin 2t)
\end{equation*}

\begin{figure}[ht]
 \centering		\includegraphics[width=.5\linewidth]{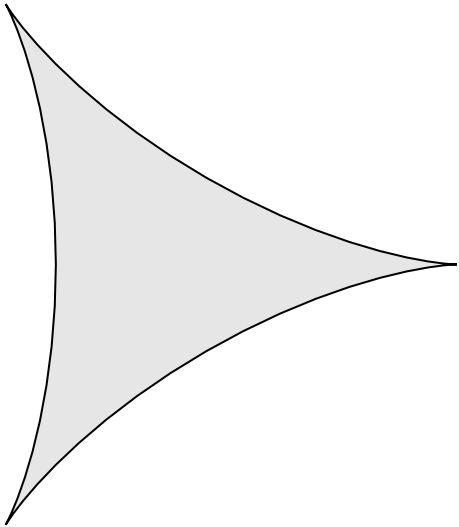}
		\caption{The deltoid domain.}
		\label{fig:Deltoide}
\end{figure}
The factor $1/3$ in the previous formulae are just here to simplify future computations, but play no fundamental rôle.
The connected component $\Omega$ of the complementary of the curve which contains $0$  is a bounded open set, that we refer to as the deltoid domain.   It turns out that there exist on this domain a one parameter family $ \LL^{(\lambda)}$ of symmetric diffusion operator which may be diagonalized in a basis of orthogonal polynomials. It was introduced in~\cite{koorn1, koorn2}, and further studied in~\cite{ zribi2013}. This is one of the 11  families of sets carrying such  diffusion operators, as described in~\cite{BOZ2013}.

In order to describe the operator, and thanks to the diffusion property~\eqref{diff.prop}, it is enough to describe $\Gamma(x,x)$, $\Gamma(x,y)$, $\Gamma(y,y)$, $ \LL^{(\lambda)}(x)$ and $ \LL^{(\lambda)}(y)$ (the  operator $\Gamma$ does not depend on $\lambda$ here).

The symmetric matrix $\bpm \Gamma(x,x)&\Gamma(x,y)\\\Gamma(y,x)&\Gamma(y,y)\epm$ is referred to in what follows as the metric associated with the operator, although properly speaking it is in fact a co-metric.
We may also  use  the complex structure of $\bbR^2\simeq \bbC$, and the complex variables $Z= x+iy$, $\bar Z= x-iy$, and it turns out that the formulas are much simpler under this description.

The operator $ \LL^{(\lambda)}$ is then described as
\beq\label{eq.deltoid}
\bcas \Gamma(Z,Z) = \bar Z-Z^2,\\ \Gamma(\bar Z, \bar Z)= Z-\bar Z^2, \\\Gamma(\bar Z, Z)= 1/2(1-Z\bar Z),\\
 \LL(Z)= -\lambda Z,  \LL(\bar Z)= -\lambda \bar Z,
\ecas
\eeq 
where $\lambda>0$ is a real parameter.

The boundary of the domain $\Omega$  turns out to the the curve with equation 
$$P(Z,\bar Z):= \Gamma(Z, \bar Z)^2-\Gamma(Z,Z)\Gamma(\bar Z, \bar Z)=0,$$and inside the domain $\Omega$, the associated metric is positive definite, so that it corresponds to some elliptic operator in $\Omega$. Moreover, for this function $P$, the boundary condition~\eqref{eq.bord} is satisfied, with
$$\Gamma(P,Z)= -3ZP, ~\Gamma(P, \bar Z)= -3\bar Z P.$$  The reversible measure associated with it, easily identified through equation~\eqref{eq.bi},  has density $C_\alpha P(Z, \bar Z)^\alpha$ with respective to the Lebesgue measure, where $\lambda = \frac{1}{2}(6\alpha+5)$, and is a probability measure exactly when $\lambda>0$ (see~\cite{zribi2013} for more details). We shall refer to this probability measure on $\Omega$  as $\mu^{(\lambda)}$.

It is quite immediate that the operator $ \LL^{(\lambda)}$ commutes with the transformation $Z\mapsto \bar Z$. 
Indeed, there is another  invariance : the transformation $Z\mapsto e^{i\theta} Z$ commutes with $ \LL^{(\lambda)}$   provided $e^{3i\theta}=1$. Therefore, everything is invariant under $Z\mapsto jZ$ and $Z\mapsto \bar j Z$, where $j$ and $\bar j$ are the third roots of unity in the complex plane.

From the form of the operator, we easily see that $ \LL^{(\lambda)}$ maps the set $\cP_n$ of polynomials with total degree $n$ (in the variables $(Z,\bar Z)$ or equivalently in the variables $(x,y)$) into itself, and being symmetric, may be diagonalized in a basis of orthogonal polynomials. We refer to~\cite{zribi2013} for a complete description of these polynomials. In what follows, and it the rest of the paper, we shall forget the dependance in $\lambda$ of the polynomials, in order to have lighter notations.

The eigenspaces of $ \LL^{(\lambda)}$  are described in~\cite{zribi2013}.  For any $(k,n)$, there is a unique polynomial $R_{n,k}$ which is a eigenvector with eigenvalue  $-\lambda_{n,k}=-((\lambda-1)(n+k) + n^2+k^2+nk)$, which is a polynomial in the variables $(Z,\bar Z)$ with real coefficients,  and has a unique highest degree  term $Z^n \bar Z^k$.   When $n\neq k$, the eigenspaces have dimension 2, and we want to distinguish between the symmetric and antisymmetric part, under the symmetry $Z\mapsto \bar Z$. We therefore chose when $n\neq k$ the basis $\frac{1}{2} (R_{n,k}+ R_{k,n}) $ and $\frac{1}{2i}(R_{n,k}-R_{k,n})$.  The following proposition summarizes  a few properties of this basis of eigenvectors.

\bprop \label{prop.defi.polyn}For any $\lambda>0$, for any $n,k\in \bbN^2$, with $k\neq n$, there are exactly, up to a sign,  two    real valued polynomials $P^{}_{n,k}(Z,\bar Z)$ and $Q^{}_{n,k}(Z, \bar Z)$, with    degree $n+k$, $P_{n,k}$ being symmetric and $Q_{n,k}$ antisymmetric under the symmetry $(Z, \bar Z)\mapsto \bar Z, Z)$,   with norm 1 in $\cL^2(\mu_\lambda)$, such that 
$$ \LL^{(\lambda)}P_{n,k}= -\lambda_{n,k} P_{n,k}, ~ \LL^{(\lambda)}Q_{n,k}= -\lambda_{n,k} Q_{n,k},$$ where 
$\lambda_{n,k} = (\lambda-1)(n+k)+ n^2+k^2+ nk$.

When $n=k$, there is exactly one such eigenvector $P_{n,n}$, with eigenvalue $-\lambda_{n,n}$, and it is symmetric in $(Z, \bar Z)$.

Moreover, $P_{n,k}$ has real coefficients, $Q_{n,k}$ purely imaginary ones, and they satisfy
\beq\label{eq.rot.j.PQ}P^{}_{n,k}(jZ, \bar j\bar Z)+ iQ_{n,k}(jZ, \bar j \bar Z)= \bar j^{n-k}\big(P^{}_{n,k}(Z, \bar Z)+ i Q_{n,k}(Z, \bar Z)\big).\eeq

\eprop

Observe that $P_{n,k}$ and $Q_{n,k}$ are real valued. We shall investigate the hypergroup property in terms of this basis.
In order to have lighter notations, we shall often write $P_{n,k}(Z)$ and $Q_{n,k}(Z)$ instead of $P_{n,k}^{}(Z, \bar Z)$ and $Q_{n,k}^{}(Z, \bar Z)$, although they are really polynomials of both variables $Z$ and $\bar Z$ (in particular, they are not harmonic in $\bbC$). Moreover, recall that we shall by convention set $Q_{n,n}=0$.

\bpf
Let $R_{n,k}(Z ,\bar Z)$ the unique eigenvector with unique highest degree term $Z^n\bar Z^k$. From the invariance under conjugacy, then  when $n\neq k$,  the conjugate of $ R_{n,k}$, that is    $\bar R_{n,k}(Z, \bar Z)$ is an eigenvector with same eigenvalue, and looking at the highest degree term, it is therefore $R_{k,n} ( Z,\bar  Z)$. Due to the conjugacy invariance, $R_{n,k}$ and $R_{k,n}$ have the same $\cL^2(\mu_\lambda)$ norm.    

For $n\neq k$, let $\hat P_{n,k}$ and $\hat Q_{n,k}$ be the symmetric  and antisymmetric eigenvectors of $ \LL^{(\lambda)}$ with dominant terms 
$\frac{1}{2}(Z^n \bar Z^k+ Z^k\bar Z^n)$ and $ \frac{-i}{2}(Z^n \bar Z^k- Z^k\bar Z^n)$ respectively, so that $\hat P_{n,k}$ and $\hat Q_{n,k}$ take real values, with $\hat Q_{n,k}$ vanishing on the real axis.  (By convention, set  $\hat Q_{n,n}=0$). The fact that $\int (P_{n,k}^{(\lambda)})^2 d\mu_\lambda=\int (P_{k,n}^{(\lambda)})^2 d\mu_\lambda$ (from the conjugacy invariance)  shows that $\hat P_{n,k}$ and $\hat Q_{n,k}$ are orthogonal. Moreover, due  to the invariance under $Z\mapsto jZ$, we know that  $\hat P_{n,k}(jZ)$ and $\hat Q_{n,k}(jZ)$ are again eigenvectors for $\LL^{(\lambda)}$ with the same eigenvalue, and therefore are linear combinations of $\hat P_{n,k}$ and $\hat Q_{n,k}$. Looking at the highest degree terms, one sees that
 \beq\label{mult.j1} \bcas \hat P_{n,k}(jZ)= (\frac{j^{n-k}+ \bar j^{n-k}}{2}) \hat P_{n,k}(Z) +i (\frac{j^{n-k}- \bar j^{n-k}}{2}) \hat Q_{n,k}(Z),\\
 \\
\hat Q_{n,k}(jZ)=-i( \frac{j^{n-k}- \bar j^{n-k}}{2}) \hat P_{n,k}(Z) + (\frac{j^{n-k}+ \bar j^{n-k}}{2})\hat  Q_{n,k}(Z).
\ecas
\eeq

In other words, $(\hat P_{n,k}+ i\hat Q_{n,k})(jZ)= \bar j^{n-k}(\hat P_{n,k}+i\hat Q_{n,k})(Z)$.
 As a consequence, since $\int \hat P_{n,k}(jZ)^2d\mu_\lambda= \int \hat P_{n,k}(Z)^2d\mu_\lambda$, one sees that $\|\hat P_{n,k}\|^2=\| \hat Q_{n,k}\|^2$  when $n-k\not\equiv 0 \mod(3)$. We do not know if this is true for $n=k \mod(3)$. Of course, the problem does not exist when $n=k$.

 We therefore may use the basis $P_{n,k}= a_{n,k} \hat P_{n,k}$ and $Q_{n,k} = b_{n,k} \hat Q_{n,k}$ as an orthonormal basis for the eigenspace associated with the eigenvalue $\lambda_{n,k}$, with $a_{n,k}= b_{n,k}$ when $n-k\not\equiv 0 \mod(3)$, and then equation~\eqref{mult.j1} translates into~\eqref{eq.rot.j.PQ}.  If $n-k\equiv 0 \mod(3)$ , then~\eqref{mult.j1} is  trivial since in this case $P_{n,k}(jZ)= P_{n,k}(Z)$, $Q_{n,k}(jZ)= Q_{n,k}(Z)$.
 \epf

There are two particular cases which are worth understanding, namely $\lambda= 1$ and $\lambda= 4$, corresponding to the parameters $\alpha= \pm1/2$.  These two models show the relation with the $A_2$ root system, and indeed our polynomials are the examples of Heckman-Opddam polynomials associated to this root system. We briefly present those two models, referring to~\cite{zribi2013} for more details.

In the first case $\lambda=1$, one sees that this operator is nothing else that the image of the Euclidean Laplace operator on $\bbR^2$ acting on the functions which are invariant under the symmetries around the lines of a regular triangular lattice. 

Indeed, consider the three unit roots of identity in $\bbC$, say $(e_1,e_2,e_3)= (1,j, \bar j)$. Then, consider the functions  $z_k: \bbC\mapsto \bbC$ which are defined  as 
\beq\label{def.zk}z_k(z) = e^{i\Re(z\bar e_k)}.
\eeq
They satisfy $|z_k|=1$ and $z_1z_2z_3=1$.

For the 2 dimensional Laplace operator, the functions $z_i$ satisfy
$$\Delta(z_i)= -z_i, \Delta(\bar z_i)= -\bar z_i,$$ and, if we denote for  $i\neq j \in \{1,2,3\}$  by  $c(i,j)$ the index in $\{1,2,3\}$ which differs from $i$ and $j$, that is   
$\{c(i,j), i,j\}= \{1,2,3\}$.
\beq\label{eq.zi.plat}\bcas \Gamma(z_i,z_i)= -z_i^2, ~\Gamma(\bar z_i,\bar z_i)= -\bar z_i^2, ~\Gamma(z_i, \bar z_i)=1\\
\\
\Gamma(z_i, z_j)= \frac{1}{2}\bar z_{c(i,j)}, ~ \Gamma(\bar z_i, \bar z_j)= \frac{1}{2} z_{c(i,j)},~
\Gamma(z_i, \bar z_j)= \frac{1}{2} z_i\bar z_j, ~i\neq j
\ecas
\eeq

Let now $Z= \frac{1}{3}(z_1+z_2+z_3)$. It is easily seen that, for the Euclidean Laplace operator on $\bbR^2$, $Z$ and $\bar Z$ satisfy the  relations~\eqref{eq.deltoid} with $\lambda=1$. Moreover, the function $Z: \bbC\mapsto \bbC$ is a diffeomorphism between the interior $\Omega_1$ of the triangle $T$ and $\Omega$, where $T$ is one of  the equilateral triangles which contains  the two edges $0$ and $4\pi/3$. The functions which are invariant under the symmetries of the triangular lattice generated by this triangle $T$ are exactly functions of $Z$. In particular, the application  $(z_1,z_2)\mapsto Z= \frac{1}{3}(z_1+z_2+ \bar z_1\bar z_2)$  maps $\cS^1\times \cS^1$  onto the closure of the deltoid domain $\Omega$. The boundary of $\Omega$ is the image under $Z$ of the boundary of the triangle, and  also of the set where 2 of the three variables $(z_1,z_2,z_3)$ coincide.  The cusps  of the deltoid model  are the images of the points $z_1=z_2=z_3$, that is $z_i= 1,j, \bar j$).

In particular, the set of all the points $Z= \frac{1}{3} (z_1+z_2+z_3)$ such that $|z_i|=1$ and $z_1z_2z_3=1$ is the closure of $\Omega$. Indeed, given $Z\in \Omega$, there exist, up to permutation,  3 unique and distinct  complex numbers $z_i$ satisfying $|z_i|=1$ and $z_1z_2z_3=1$ such that $Z= \frac{1}{3}(z_1+z_2+z_3)$. They are the three distinct  roots  of the equation  $P(X)= X^3-3ZX^2+ 3\bar ZX-1 =0$. It is worth for this  to observe that the discriminant of $P$ is,  up to a numerical constant, equal to $P(Z,\bar Z)$, and therefore does not vanish in $\Omega$. 

The second description comes  from the Casimir operator on $SU(3)$. 
For $SU(d)$, one may describe this operator through its action on the various entries $z_{ij}$ of the matrices $g\in SU(d)$. Up  to some scaling factor, one has
$$\bcas  \Gamma^{SU(d)}(z_{ij}, z_{kl})= z_{ij}z_{kl}-dz_{il}z_{kj}, \\
 \Gamma^{SU(d)}(\bar z_{ij},\bar  z_{kl})= \bar z_{ij}\bar z_{kl}-d\bar z_{il}\bar z_{kj},\\
 \Gamma^{SU(d)}(z_{ij}, \bar z_{kl})= d\delta_{ik}\delta_{jl}-z_{ij}\bar z_{kl}
 \\  \LL^{SU(d)}(z_{ij})= -(d^2-1) z_{ij}, ~ \LL^{SU(d)}(\bar z_{ij})= -(d^2-1) \bar z_{ij}\ecas
 $$
Then, choosing $Z= \frac{1}{3}(z_{11}+z_{22}+z_{33})$, one obtains the relations~\eqref{eq.deltoid} with $\lambda= 4$ when applying  $\frac{3}{4} \LL^{SU(3)}$ to functions of $Z$ and $\bar Z$.
For this particular case, we may make full use of the group structure of $SU(3)$ to obtain the hypergroup property for the deltoid model, at any of the points $Z=1, Z= e^{2i\pi/3}$, $Z= e^{4i\pi/3}$. It is enough to follow the scheme described in Section~\ref{sec.hyperg}, chosing $Y= SU(3)$, $\pi$ being the application which to some matrix $g\in SU(3)$ associates $\frac{1}{3}\tr(g)$, and $\Phi$ being $g\mapsto g_0g$, for any $g_0\in SU(3)$. . If $\pi(g)=1$ for example, then $g= \Id$, and the conditional of $\pi(\Phi(g))$ knowing that $\pi(g)=1$ is the Dirac mass at $\pi(g_0)$. However, the eigenspaces being two dimensional, this introduces some extra complexity that we shall examine in Section~\ref{sec.hgp.deltoid1}.

Similarly, using the representation for $\lambda= 1$, we may also prove the hypergroup property in this case, using the translations in $\bbR^2$. 

Our aim in what follows is to propose another $6$-dimensional model which projects onto the deltoid model in the general case (that is for $\lambda\neq 1,4$), and on which we have enough symmetry to use the machinery described in Section~\ref{sec.hyperg}.

 \section{A $6$-dimensional model for the deltoid\label{6.dim.model}}
 
 In this section, we construct a symmetric diffusion operator $\hat  \LL_\lambda$ in dimension $6$ (or more precisely on an open bounded  set $\Omega_1\subset\bbC^3$), such that the image  of the operator $\hat  \LL_{\lambda}$ under the projection $\pi: \bbC^3\mapsto \bbC$  which is $\pi(z_1,z_2,z_3)= \frac{1}{3}(z_1+z_2+z_3)$ is exactly our deltoid model with parameter $\lambda$. Moreover, this operator may be diagonalized in a complete system of orthogonal polynomials for its reversible measure.

We  consider a diffusion operator in $\bbC^3$, acting on  $3$ complex variables $z_1,z_2,z_3$,  defined as follows. 
We set

\beq\label{defL}\bcas
\hat\Gamma(z_i,z_j)=\frac{3}{2} \bar z_{c(i,j)}-z_iz_j & \hbox{ for $i\neq j$}\\
\hat\Gamma(\bar z_i,\bar z_j)=\frac{3}{2}  z_{c(i,j)}-\bar z_i\bar z_j & \hbox{ for $i\neq j$}\\
\hat\Gamma(z_i,z_i)= -z_i^2&\\
\hat\Gamma(\bar z_i,\bar z_i)= -\bar z_i^2&\\
\hat\Gamma(z_i, \bar z_j)= \frac{3}{2}\delta_{ij}-\frac{1}{2}z_i\bar z_j\\
\hat \LL(z_i)= -\lambda z_i, ~\hat \LL(\bar z_i)= -\lambda \bar z_i,
\ecas
\eeq
where as before $c(i,j)$ the index in $\{1,2,3\}$ which differs from $i$ and $j$.
It is worth to observe that those relations fit with~\eqref{eq.zi.plat} whenever $|z_i|= 1$ and $z_1z_2z_3=1$, so that those equations may be seen as an extension of~\eqref{eq.zi.plat} in the interior of some domain bounded by those equations. We shall see hat it is indeed the case, when we shall describe the domain.

The first task is to observe that $\hat \LL$  projects onto the deltoid model in the sense of Section~\ref{sec.hyperg}, with the same parameter $\lambda$.

\bprop For $z= (z_1,z_2,z_3)\in \bbC^3$, let $\pi(z):=Z= \frac{1}{3}(z_1+z_2+z_3)$. Then, one has
$$\bcas \hat \Gamma(Z, Z)= \bar Z- Z^2\\
\hat \Gamma(\bar Z,\bar Z)= Z-\bar Z^2\\
 \hat \Gamma(Z, \bar Z= \frac{1}{2}\big(1-Z\bar Z)\big)\\
 \hat  \LL(Z)= -\lambda Z,  ~\hat  \LL (\bar Z)= -\lambda \bar Z
 \ecas
$$

\eprop
This comes from an immediate computation. We see that these computations fit with formulae~\eqref{eq.deltoid}.

The next  task is to show that this corresponds to some elliptic symmetric diffusion generator in some bounded domain $\Omega_1\subset \bbC^3$.

For this, let us introduce the following notations. We write $z_j= r_je^{i\theta_j}$, with $r_j=|z_j|$, and 
$$\bcas \theta= \theta_1+\theta_2+ \theta_3\\
S_1= r_1^2+r_2^2+r_3^3,\\ S_2= r_1^4+ r_2^4+ r_3^4\\ 
\sigma= r_1r_2r_3\cos(\theta)
\ecas
$$ and let 
\beq\label{def.P1.P2}\bcas
 P_1(z_1,z_2,z_3)&= 2-(S_1+1)^2 + 2S_2+ 8\sigma.\\
P_2(z_1,z_2,z_3)&= 2(S_2-1)- (S_1-1)^2.
\ecas
\eeq

Let us describe the domain first. Let $D$ be the determinant of the matrix $\hat\Gamma$ (in both coordinates $(z_i, \bar z_i)$). Then, define $\Omega_1$ to be  the connected component of $\bbC^3\setminus \{D=0\}$ which contains $(0,0,0)$. 

\bprop We have
\benum
\item $D= \frac{3^5}{2^6}P_1P_2$;
\item  The operator   $\hat \LL$ is elliptic  in $\Omega_1$;
\item $\Omega_1$ is bounded;
\item $P_1>0$ on $\Omega_1$ and $P_2<0$ on $\Omega_1$.
\item The boundary of $\Omega_1$ is included in $\{P_1=0\}$.

\eenum
\eprop
\bpf The first item results from a direct computation. Since it is  straightforward although quite technical, we recommend that the reader uses a computer program to check it.

For the second item, let us observe that the metric at $z= (0,0,0)$ is the classical Euclidean metric of $\bbC^3$, up to some scaling factor. Therefore, it  remains elliptic as long as no eigenvalue of the metric vanishes,  that is in  the  connected component of the set $\{D\neq 0\}$ which contains $\{0,0,0\}$. It is  therefore elliptic in $\Omega_1$.

For the third assertion, observe that as long as $\hat \Gamma$ is positive definite, then one has 
 $\Gamma(z_i, \bar z_i)^2> \Gamma(z_i,z_i)\Gamma(\bar z_i, \bar z_i)$.  This translates into 
$(3+r_i^2)(1-r_i^2)>0$ so that  $r_i\leq 1$ in $\Omega_1$.

 We have $P_2(0,0,0)<0$ and $P_1(0,0,0) >0$, so that, in $\Omega_1$,  $P_1>0$ and $P_2<0$.  
 
 It remains to prove  the last assertion.  Any point  $z\in \partial\Omega_1$ satisfies $P_1(z)P_2(z)=0$. Let $z= (z_1,z_2,z_3)\in \partial\Omega_1$ such that $P_2(z)=0$.  We shall prove that we also have $P_1(z)=0$ at this point. 
 
  We know that $S_1\geq 0$, and it is immediate that   $\frac{1}{3}S_1^2\leq S_2\leq S_1^2$. From the upper bound, we get $P_2\leq ( S_1-1)(S_1+3)$. Therefore,    at a point where $P_2=0$, one has $S_1\geq 1$. At such a point, we also have 
 $P_1= 4(1-S_1+2\sigma)$. 
 
 We now want to estimate the maximum possible value of $r_1r_2r_3$ on the set where $P_2=0$ and $0\leq r_i\leq 1$. We shall see that this maximum value is $(S_1-1)/2$.  Setting $x_i= r_i^2$, it amounts to look for the maximum value of $x_1x_2x_3$ on the set where $S_2:= \sum_i x_i^2= 1+ \frac{1}{2}(S-1)^2$, where $S= \sum_i x_i$ and $0\leq x_i\leq 1$.  Looking at the Lagrange multiplier, one sees that, if the maximum is not attained at the boundary of the set, then one must have 
 $$\bcas x_2x_3= \lambda(1+x_1-x_2-x_3)\\
 x_1x_3= \lambda(1+x_2-x_1-x_3)\\
 x_1x_2= \lambda(1+x_3-x_1-x_2)\ecas
 $$ from which using $S_2= 1+ \frac{1}{2}(S-1)^2$, one deduces that 
 $3\lambda(3-S)= (S-1)(S+3)$ and therefore $\lambda>0$. Multiplying the first equation by $x_1$ we find then that 
 $$\pi:= x_1x_2x_3= \lambda(x_1^2+x_1+2\lambda(x_1-1)).$$
 The same relation being true for $x_2$ and $x_3$,  the values for $x_i$ must lie among the two solutions of the equation
 $$\pi=  \lambda(x^2+x+2\lambda(x-1)).$$  But $\pi>0$ and $\lambda>0$.  Unless they are all equal, they may not be all positive. Therefore, either the maximum is attained at a point $(x,x,x)$, in which case this common value is $1$ and  $(S_1-1)/2$, either one of the values for $x_i$ is $1$, and the fact that $P_2=0$ implies that the two other values are equal, in which case again the value is $(S-1)/2$. 
 
Therefore, when $P_2=0$,  $P_1\leq 1-S_1+ \frac{1}{2}(S_1-1)= \frac{1}{2}(1-S_1)\leq 0$. 
 Since $P_1>0$ in $\Omega_1$, then any point in $\partial \Omega_1$ which satisfies $P_2(z)=0$ also satisfies $P_1(z)=0$. 

\epf 

\brmq\label{rmq.decomp.P1} It may be worth to observe that $P_1$ may take a somewhat simpler form. Setting 
$$\bcas \sigma_0= r_1+r_2+r_3\\
\sigma_1= -r_1+r_2+r_3\\
\sigma_2= r_1-r_2+r_3\\
\sigma_3= r_1+r_2-r_3\ecas
$$ and 
\beq\label{form.A1.B1}\bcas
S= (1+\sigma_0)(1-\sigma_1)(1-\sigma_2)(1-\sigma_3)\\
D= (1-\sigma_0)(1+\sigma_1)(1+\sigma_2)(1+\sigma_3)
\ecas
\eeq
one may write 
$$P_1= S\cos^2(\theta/2)+ D\sin^2(\theta/2).$$

\ermq

Our next result shows that indeed the operator $\hat \LL$ defined on $\Omega_1$ satisfies the boundary condition~\eqref{eq.bord}.

\bprop~\label{prop.L.measure.bord}
\benum
\item  With $P_1$ defined in~\eqref{def.P1.P2}, one has 
\beq \label{cond.poly.delt2}\hat\Gamma(\log(P_1), z_i)= -3z_i, \hat\Gamma(\log(P_1), \bar z_i)= -3\bar z_i.\eeq

\item For $i= 1, \cdots, 3$,  one has 
\beq  \bcas
\sum_{j=1}^3 \partial_{z_j} \hat\Gamma(z_j,z_i)+ \partial_{\bar z_j}\hat \Gamma(\bar z_j,z_i) = -\frac{11}{2} z_i\\
\sum_{j=1}^3 \partial_{z_j}\hat \Gamma(z_j,\bar z_i)+ \partial_{\bar z_j} \hat\Gamma(\bar z_j,\bar z_i) = -\frac{11}{2} \bar z_i
\ecas
\eeq

\eenum
\eprop

\bpf These   formulas may be checked  with a direct and tedious  computation.  However, we have no simple interpretation, beyond mere calculus, why this is true. The first one is the condition required for the operator $\hat  \LL$ on $\Omega_1$ to be diagonalizable in a system of orthogonal polynomials. The second will be used to identify the reversible measure.

\epf

\brmq Together with the fact that $\Gamma(z_i, z_j)$, $\Gamma(z_i, \bar z_j)$ and $\Gamma(\bar z_i, \bar z_j)$ are polynomials in the variables $(z_i,\bar z_i)$, and following~\cite{BOZ2013},  equation~\eqref{cond.poly.delt2} is the key identity which insures that  the operator $\hat \LL$ may be be diagonalized in a basis of orthogonal polynomials.

\ermq
We now prove the following lemma 
\blem \label{lem.integr}For any $\beta>-1$, $\int_{\Omega_1} P_1^{\beta} dx<\infty$, where $dx$ is the Lebesgue measure on $\Omega_1\subset \bbC^3$.

\elem
Before proving Lemma~\ref{lem.integr}, we deduce the  following
\bcor\label{cor.inv.meas} For any $\lambda>5/2$, the operator $\hat \LL^{(\lambda)}$ is reversible with respect to the probability  measure $C_\beta P_1^\beta dx$, where $\beta = \frac{1}{6}(2\lambda-11)$, $dx$ is the Lebesgue measure in $\Omega_1$,  and $C_\beta$ is the normalizing constant.

\ecor 

\bpf(Of Corollary~\ref{cor.inv.meas})
This is a direct consequence of the general formula~\eqref{eq.bi} and Proposition~\ref{prop.L.measure.bord}.

\epf

\brmq 
It is worth to observe that $\beta = -1$ corresponds to $\lambda= 5/2$, while 
$$\hat \LL(P_1)=-4(\lambda+2)P_1+ (2\lambda-5)(P_2-3P_1)/6.$$

This suggest that, for this precise value $\lambda= 5/2$, the measure is concentrated on the surface $\{P_1=0\}$, and the associated process  lives indeed on the boundary of the domain $\Omega_1$.

Moreover, this limit $\beta = 5/2$ corresponds for the projected model on the deltoid domain to the case $\alpha=0$, that is when the reversible measure of the image operator on $\Omega$ is the Lebesgue measure.

\ermq

 \bpf (Of Lemma~\ref{lem.integr})
 
 When looking at the  behavior of $P_1$ near a regular point of the boundary $\{P_1=0\}$, it is clear that the condition $\beta >-1$ is necessary for $P_1^\beta $ to be locally integrable near such a point, in the domain $\{P_1>0\}$. But we have also to consider the behavior near  singular points. 
Set as before 
 $$\bcas S_1= r_1^2+r_2^2+r_3^2,\\ S_2= r_1^4+r_2^4+r_3^4, \\\tau = r_1r_2r_3. \ecas$$
and, writing  $z_j= r_je^{i\theta_j}$, set $t= \theta_1+\theta_2+\theta_3$. Then we have   

 $$P_1= 2-(S_1+1) ^2+ 2S_2+ 8\tau \cos(t),$$ and we want to show that for $\beta>-1$, $P_1^\beta$ is locally integrable  with respect to the measure $r_1r_2r_3dr_1dr_2dr_3dt$.
 on the domain 
 $0<r_i<1$, $t\in (0,2\pi)$,  and $P_1>0$. 
 We set 
 $$A=  2-(S_1+1) ^2+ 2S_2, ~B= 8\tau, $$
 The first thing is to compute  
 \beqna I(A,B)= &&\int_{A+B\cos(t)>0, t\in (0,2\pi)} (A+B\cos(t))^{\beta } dt\\\label{eq.int.P1}
 &&=2\int_{-1}^1\I_{A+Bu>0} (A+Bu)^{\beta } \frac{du}{\sqrt{1-u^2}},
 \eeqna
 that we want to estimate up to some constants depending only on $\beta $. For this we write   $I\simeq J$ when $c_\beta \leq \frac{I}{J} \leq C_\beta $, for two positive constants $c_\beta$  and $C_\beta $.
 
We  cut the integral in~\eqref{eq.int.P1}  into $\int_0^1$ and $\int_{-1}^0$. Then we  change variables to write both integrals  as $\int_0^1 $, and again  change variables $v\mapsto 1-v$. Finally,  we get
$$I(A,B)\simeq \int_0^1\I_{A+B-Bv>0}(A+B-Bv)^{\beta }\frac{dv}{v^{1/2}}+ \int_0^1\I_{A-B+Bv>0}(A-B+Bv)^{\beta }\frac{dv}{v^{1/2}}.$$
Write this $I(A,B)= I_1+I_2$.

 Then
 $$\bcas I_1(A,B)= 0\hbox{ if } A+B<0, \\
  I_1(A,B)\simeq (A+B)^{\beta +1/2}B^{-1/2}~\hbox{ if } A+B>0 
 \ecas 
 $$
 and 
 $$\bcas 
 I_2=0 \hbox{ if } A<0, \\ 
 I_2\simeq A^{\beta +1}B^{-1} \hbox{ if } 0<A<B\\ 
I_2\simeq A^{\beta } \hbox{ if } 0<B<A\ecas
$$

 We now have to consider the integral of $I_1$ and $I_2$ on the image of the domain $\Omega_1$ through the projection $(z_1,z_2,z_3)\mapsto (r_1,r_2,r_3)$,  that is the integral of $I_1$ and $I_2$  with respect to the measure $r_1r_2r_3 \,dr_1dr_2dr_3$, on the domain $\{0\leq r_i\leq 1\} \cap \{A+B>0\}\cap \{A-B>0\}$.  An easy inspection show that $A+B=S$ and $A-B= D$, where $S$ and $D$ are given in formulae~\eqref{form.A1.B1}.
We see that $I_1\geq C( l_1l_2l_3)^{\beta +1/2}B^{-1/2}$, where $l_i $ are 3 independent affine forms, restricted to a domain where they are positive, and $C$ is a constant. Such a function is integrable with respect to the measure $B^{1/2}\un_{0<r_1\leq 1}\un_{0<r_2\leq 1}\un_{0<r_3\leq 1}dr_1dr_2dr_3$. This integral is finite when $\beta >-1$.

The only problem comes from $\int_{ 0<B<A, 0<r_i\leq 1}  A^{\beta } r_1r_2r_3 dr_1dr_2dr_3$. We see that one may reduce to a neighborhood of a point where both $A$ and $B$ vanish, which means that (up to a permutation of the indices) $r_1=0$, $r_2+r_3=1$. Once again, integrating in a area where $0<< r_2<<1$ causes no problem, and it remains to consider the integral on an area where $r_2=0, r_3=1$ ( or symmetrically). Then, in a neighborhood of this point,  we may bound $A\geq A-B= D$, and consider the integral of $D^\beta $ on $D>0$. In this neighborhood, with the notations of Remark~\ref{rmq.decomp.P1}, $D$ is equivalent to $(1-\sigma_0)(1+\sigma_3)$, which is the product of two independent linear forms (in the variables $r_1,r_2,r_3$) , and is locally integrable as soon as $\beta >-1$.

\epf

 \section{ Hypergroup property for the deltoid model\label{sec.hgp.deltoid1}}
 
 For $\lambda>5/2$, the operator $\hat{ \LL}^{(\lambda)}$ provides a symmetric diffusion operator on the set $\Omega_1$ which projects onto the deltoid model $ \LL^{(\lambda)}$ with the same parameter through the map $\pi(z_1,z_2,z_3)=\frac{1}{3}( z_1+z_2+z_3)$. 
 
  We now  observe that, for any $\theta= (\theta_1, \theta_2)\in \bbR^2$,  the operator $\hat \LL$ is invariant under the transformations $\Phi_{\theta}(z_1,z_2,z_3)= (e^{i\theta_1} z_1, e^{i\theta_2} z_2, e^{-i(\theta_1+\theta_2)} z_3)$.  It is worth to notice that $\Phi_{\theta_1}\Phi_{\theta_2}= \Phi_{\theta_1+\theta_2}$, and that the measure $\mu_\alpha= P_1^\alpha dz$ is invariant under $\Phi_\theta$.  Therefore, if we define $\Phi_\theta(f)(z)= f\big(\Phi_\theta(z)\big)$, one see that the adjoint of $\Phi_\theta$ is $\Phi_{-\theta}$
 
 To fit with the situation described in Section~\ref{sec.hyperg}, it remains to identify the point $x_0$ of Proposition~\ref{prop.hpgp.gal}. Let then $x_0= \pi(1,1,1)= 1$: it belongs to the boundary of $\Omega$. According to the analysis of the flat model ($\lambda=1$), $\theta\mapsto \Phi_{\theta}(1,1,1)$ is onto the deltoid domain $\Omega$, and, when $\theta$ varies in $\bbR^2$.  Indeed, we have seen that any point in $\Omega$ may be written as $\frac{1}{3}(z_1+z_2+z_3)$, where $|z_i|=1$ and $z_1z_2z_3=1$, which corresponds to points $\Phi_\theta(1,1,1)$.
 
 Since $|z_i|^2\leq 1$ for any point in $\Omega$, then, $\forall z\in \Omega$,  $|\pi(z)|\leq 1$,  and if $\pi(z)=1$, then $z_1=z_1=z_3= 1$. Therefore, $\pi(\Phi_{\theta}(z_1,z_2,z_3))= \Phi_{\theta}(1,1,1)$. This shows that the conditional law of $\pi\big(\Phi_\theta(z)\big)$ when $\pi(z)= 1$ is a Dirac mass at 
$\Phi_{\theta}(1,1,1)$.  

  Indeed, we may as well chose  for $x_0$ the image of the points $(j,j,j)$ or $(\bar j, \bar j, \bar j)$, and those three points correspond to the cusps of the deltoid  curve.
Define $Z(\theta)= \pi((\Phi_\theta(1,1,1))$. Then,   $\pi(\Phi_\theta (j,j,j))= jZ(\theta)$, and   $\pi(\Phi_\theta (\bar j,\bar j, \bar j ))= \bar jZ(\theta)$.  Moreover, the conditional law of  $z= \frac{1}{3}(z_1+z_2+z_3)$ knowing that $\pi(z)= 1,j,\bar j$ is a Dirac mass at $Z(\theta)$, $jZ(\theta)$ and $\bar j Z(\theta)$ respectively.  The point $Z(\theta)$ is in the interior of the deltoid domain as soon as $\theta_1\neq \theta_2 \mod(2\pi)$ and $2\theta_1\neq -\theta_2 \mod(2\pi)$. Observe also that $Z(-\theta)= \bar Z(\theta)$.

To apply the method described in Section~\ref{sec.hyperg}, the only problem is that the eigenspaces associated with $ \LL^{(\lambda)}$ have dimension 2, when $n\neq k$. They are invariant under the symmetry $Z\mapsto \bar Z$, 
and the model $\hat  \LL^{(\lambda)}$ also shares the symmetry $\cS : (z_1,z_2,z_3)\mapsto (\bar z_1,\bar z_2, \bar z_3)$. So instead of looking at eigenvectors of $ \LL^{(\lambda)}$ alone, one may look at eigenvectors of $ \LL^{(\lambda)}$ which are symmetric, or antisymmetric, through the transformation $Z\mapsto \bar Z$, which leads us to consider the basis $(P_{n,k}, Q_{n,k})$ described in section~\ref{sec.delt.mod}.

The transformation $z\mapsto \Phi_\theta(z)$ is not invariant under  the symmetry $\cS$. We have $\overline{\Phi_\theta(z)}= \Phi_{-\theta} (\bar z)$.
Observe also that the operators $\Phi_\theta(f)(z)= f(\Phi_\theta (z)) $ satisfy $\lag \Phi_\theta f, g\rag= \lag f, \Phi_{-\theta } g \rag$, which comes from the fact that the measure is invariant under the symmetry~$\cS$.

From the general scheme described in Section~\ref{sec.hyperg}, one sees that, writing the Markov operator  $K_\theta(f)= \bbE( f(\pi(\Phi_\theta z))/ \pi(z)= Z)$, that both $  K_\theta(P_{n,k})$ and  $K_\theta(Q_{n,k})$ belong to the eigenspace associated to the eigenvalue $\lambda_{n,k}$, defined in Proposition~\ref{prop.defi.polyn}. We have

\bprop With $K_\theta(f)(Z)= \bbE( f(\pi(\Phi_\theta z))/ \pi(z)= Z)$, one has, for any $\lambda\geq 5/2$, 
$$ \bpm K_\theta(P_{n,k})\\K_\theta(Q_{n,k})\epm = \bpm \alpha_{n,k}(\theta)  & \beta_{n,k}(\theta)\\ \gamma_{n,k}(\theta)& \delta_{n,k}(\theta)\epm \bpm P_{n,k}\\Q_{n,k}\epm,$$

with 
$$\bcas \alpha_{n,k}(\theta)=   \frac{P_{n,k}(Z(\theta))}{P_{n,k}(1)};\\
\\
\gamma_{n,k}(\theta)= -\beta_{n,k}(\theta) =\frac{Q_{n,k}(Z(\theta))}{P_{n,k}(1)}.
\ecas
$$ 
Moreover, one has 
$\alpha_{n,k}(\theta)= \alpha_{n,k}(-\theta)$, $\delta_{n,k}(-\theta)= \delta_{n,k}(\theta)$,  and $\beta_{n,k}(\theta)= -\beta_{n,k}(-\theta)$.
\eprop
\bpf It is enough to check this property for $\lambda>5/2$, since everything is continuous and we would have the same property in the limit $\lambda= 5/2$.

Equation $\lag K_\theta(f), g\rag = \lag f, K_{-\theta}(  g)\rag$ proves the parity relations, choosing $f$ and $g$ either $P_{n,k}$ or $Q_{n,k}$.

If we apply $K_\theta (f)(1)= f\big(Z(\theta)\big)$,  and using the fact that $Q_{n,k}(1)=0$, we get $P_{n,k}(1)\alpha_{n,k}(\theta)= P_{n,k}\big(Z(\theta)\big)$, and $P_{n,k}(1)\gamma_{n,k}(\theta)= Q_{n,k}\big(Z(\theta)\big)$.  This shows that $P_{n,k}(1)$ may not vanish (since it would imply that $P_{n,k}=0$ everywhere), and leads to the representation formula for $\alpha_{n,k}$ and $\gamma_{n,k}$. The formula for $\beta_{n,k}$ follows by symmetry, since $Z(-\theta)= \bar  Z(\theta)$ and $Q_{n,k}(\bar Z)= -Q_{n,k}(Z)$.

\epf

\brmq The value of $\delta_{n,k}(\theta)$ however is more difficult to obtain.  One may apply the identity $K_\theta(f)(j)= f(jZ(\theta))$ and formulae~\eqref{mult.j1} whenever $n-k\neq 0 \mod(3)$, to get
$$\delta_{n,k}(\theta)= \frac{\cos(2(n-k)\pi/3)}{\sin(2(n-k)\pi/3)} \frac{P_{n,k}(Z(\theta))}{P_{n,k}(1)}.$$ This method, however, does not provide any information on $\delta_{n,k}(\theta)$ when $n-k \equiv 0 \mod(3)$.

\ermq

The point $x_0=1$ corresponds to one of the 3 cusps of the deltoid curve, and is the image of one of the points of $\cS^1\times \cS^1\times \cS^1$ where $z_1=z_2=z_3$ and $z_1z_2z_3=1$. There are 3 such points, corresponding to the three cusps of the deltoid, and we could have similarly proved  the hypergroup property for any of those points, but for another polynomial basis.  The choice of the point $x_0=1$ corresponds to the choice of the basis  $(P_{n,k}, Q_{n,k})$ such that, under the symmetry  $\cS :Z\mapsto \bar Z$,  $SP_{n,k}= P_{n,k}$ and $SQ_{n,k}= -Q_{n,k}$.  

We may as well consider symmetries which leave the two other cusps invariant, and this provides new bases for the eigenspace in which the operator $K_\theta$ has a similar expression.
A polynomial $R$  is symmetric with respect to the symmetry through the $j$ axis if $R(\bar j \bar Z)= R(Z)$.  Then, thanks to equation~\eqref{eq.rot.j.PQ}, a basis $(R_{n,k}, S_{n,k})$ of the eigenspace associated with $\lambda_{n,k}$ for which the first element is symmetric and the second antisymmetric under the symmetry around the $j$ axis would be 
$R_{n,k}+ iS_{n,k}= j^{n-k}(P_{n,k}+iQ_{n,k})$.
One also have the hypergroup property for the family $(R_{n,k}, S_{n,k})$ and similarly for the family $(U_{n,k},V_{n,k})$ corresponding to the point $\bar j$.This leads to other representations of the Markov operator $K_\theta$.

  As usual, the operators $K_\theta$ lead to a full representation of Markov kernels.

\bthm  Let $K$ be a symmetric  Markov operator, bounded in $\cL^2(\mu^{(\lambda)})$, with $\lambda \geq 5/2$. Assume that $K$  commutes with $ \LL^{(\lambda)}$. Then, with the notations of Proposition~\ref{prop.defi.polyn}, for any $n,k\in \bbN^2$, it satisfies 
\beq\label{eq.rep.K} \bpm K(P_{n,k})\\ K(Q_{n,k})\epm = \bpm a_{n,k}& b_{n,k}\\b_{n,k} &c_{n,k}\epm \bpm P_{n,k}\\Q_{n,k}\epm\eeq
and there exists a probability measure $\nu_1$ on the deltoid domain $\Omega$ such that
 \beq\label{repr.1}a_{n,k} = \int \frac{P_{n,k}(z)}{P_{n,k}(1)} \nu(dz), ~b_{n,k} = \int \frac{Q_{n,k}(z)}{P_{n,k}(1)} \nu(dz).\eeq 
 \ethm.
 
 \bpf

We follow the lines of the proof described in Section~\ref{sec.hyperg}. Equation~\eqref{eq.rep.K} is immediate from the fact that $K$ commute with $ \LL^{(\lambda)}$ and the description of the eigenspaces of $ \LL^{(\lambda)}$ given in Proposition~\ref{prop.defi.polyn}. Extending the operator $K$ to act on probability measures, we choose   $\nu_1(dz)= K(\delta_1)$.  If $P_t$ denotes the heat kernel associated with $ \LL^{(\lambda)}$, let $\nu_t= K\big(P_t(\delta_1)\big)$.  Following~\cite{BakryZribiCDdeltoide}, we know that $P_t(\delta_1 )$ has a bounded density with respect to $\mu^{(\lambda)}$, which may be written as 
$\sum_{n,k} e^{-\lambda_{n,k}t} P_{n,k}(z)P_{n,k}(1)$, where this simplified form comes from the fact that $Q_{n,k} (1)=0$. Then, the density  $h_t$  of $\nu_t$  with respect to $\mu^{(\lambda)}$ may be written as 
 $$h_t=\sum_{n,k} e^{-\lambda_{n,k} t}P_{n,k}(1) (a_{n,k} P_{n,k}(Z)+ b_{n,k} Q_{n,k}(Z)),$$ and we see that
 $$\int \frac{P_{n,k}(Z)}{P_{n,k}(1)} d\nu_t= e^{-\lambda_{n,k}t} a_{n,k}, ~\int \frac{Q_{n,k}(Z)}{P_{n,k}(1)} d\nu_t=e^{-\lambda_{n,k}t} b_{n,k}.$$
 Since $\nu_t$ converges to $\nu= K(\delta_1)$ when $t\to 0$,  we get the result in the limit.
 \epf

 The previous representation relies in an essential way on the fact that $Q_{n,k}(1)=0$. This choice comes  from the symmetry properties of the operator under $Z\mapsto \bar Z$, that is the symmetry around the real axis. 
 
 In order to get informations about the coefficients $c_{n,k}$, one may use the invariance of the model under $Z\mapsto jZ$ and $Z\mapsto \bar j Z$ whenever $n-k\not\equiv 0 \mod (3)$. One may use similarly the invariance through multiplication by $j$ and $\bar j$ and symmetries with the corresponding axis.
 This means that a similar presentation is valid in the two other bases $(R_{n,k}, S_{n,k})$  and $(U_{n,k}, V_{n,k})$ (using the symmetries leaving $j$ and $\bar j$ invariant respectively).
 
 In this new basis, the matrix of the operator is unchanged when $n-k\equiv 0 \mod(3)$, while  the matrix of the operator $K$ in this new basis becomes,  
 $$\frac{1}{4} \bpm a_{n,k} +  2\epsilon \sqrt{3} b_{n,k}+ 3c_{n,k}& -\epsilon \sqrt{3} a_{n,k} - 2b_{n,k} + \epsilon \sqrt{3} c_{n,k}\\
-\epsilon \sqrt{3}a_{n,k} - 2b_{n,k} +\epsilon  \sqrt{3}c_{n,k}&3a_{n,k}-  2\epsilon \sqrt{3} b_{n,k}+ c_{n,k}\epm
$$
where $\eps= 1$ when $n-k\equiv 1 \mod(3)$ and $\eps= -1$ when $n-k\equiv 2 \mod(3)$.

If we observe that $(R_{n,k}+ iS_{n,k}) (jZ)= (P_{n,k}+ iQ_{n,k})(Z)$, so that $R_{n,k}(j)= U_{n,k}(\bar j)= P_{n,k}(1)$, we get a new representation, with the measure $\nu_1= K(\delta_j)$, when $n-k\not\equiv 0\mod(3)$
$$\frac{1}{4}(a_{n,k}+2\epsilon \sqrt{3} b_{n,k} + 3 c_{n,k}) = \int \frac{R_{n,k}(z)}{P_{n,k}(1)} \nu_1(dz)= \int \frac{U_{n,k}(z)}{P_{n,k}(1)} \nu_2(dz) $$
and $$\frac{1}{4}( -\epsilon \sqrt{3} a_{n,k} - 2b_{n,k} + \epsilon \sqrt{3} c_{n,k})= \int \frac{S_{n,k}(z)}{P_{n,k}(1)} \nu_1(dz)= \int \frac{V_{n,k}(z)}{P_{n,k}(1)} \nu_2(dz) .$$

This may be rewritten as 
$$\bcas\dis -\frac{a_{n,k}}{2}-\eps \frac{\sqrt{3}}{2}b_{n,k}= \int \frac{P_{n,k}(z)}{P_{n,k}(1)} d\nu_1(z)\\
\\
\dis -\frac{b_{n,k}}{2} -\eps\frac{\sqrt{3}}{2} c_{n,k} = \int \frac{Q_{n,k}(z)}{P_{n,k}(1)} d\nu_1(z).\ecas$$

 Changing $j$ into $\bar j$ amounts to change $\epsilon$ into $-\epsilon$ in the previous formulas, with the measure $\nu_2= K(\delta_{\bar j})$.
 
 This in turn provides a representation of $c_{n,k}$ when $n-k\not\equiv k \mod(3)$, of the form
 $$c_{n,k}=\eps\int\frac{Q_{n,k}(z)}{P_{n,k}(1)} \Big( \frac{2}{\sqrt{3}}d\nu_1(z)-\frac{1}{\sqrt{3}}d\nu(z)\Big).$$ 
 
 Unfortunately, this does not carry any information about $c_{n,k}$ when $n\equiv k \mod(3)$.
 
\brmq As we already observed,  for  $\lambda= 5/2$, the $6$ dimensional model is in fact carried by the algebraic hypersurface $\{P_1=0\}$, and in fact is a $5$-dimensional model. 
However, for other values of $\lambda$ (with the sole exception of $\lambda=1$), we do not know if the property is true. It would be intersecting to construct lower dimensional models for these values, but this seems quite hard.
\ermq

\brmq It would be interesting to have an explicit expression for the kernel $K_\theta(x,dy)$, in order to have explicit representation for the product formula~\eqref{eq.product.formula}. Unfortunately, the law of $(Z, R_\theta(Z))$ is   already apparently quite out of reach through simple formulas. 

\ermq
\section{ Projection of the deltoid model and the $G_2$-root system\label{sec.G2}}

As we saw in Section~\ref{sec.hgp.deltoid1},  the fact that the eigenspaces for $ \LL^{(\lambda)}$ are two dimensional introduce extra complexity in the representation of the eigenvalues of  the Markov operators which commute with $ \LL^{(\lambda)}$. This works much better if we concentrate on functions which are symmetric in $(Z, \bar Z)$, which correspond to the symmetric polynomials $P_{n,k}(Z, \bar Z)$. . It turns out that these symmetric polynomials are again orthogonal polynomials corresponding to another bounded set $\Omega_2\subset \bbR^2$, bounded by a cuspidal cubic and a parabola, tangent to each other at the second order. Going back to the triangle model, remember that  the deltoid model in the case $\lambda=1$ corresponds  to functions which are invariant under the symmetries of a triangular lattice, corresponding to the root system $A_2$. Adding this new invariance under $Z\mapsto \bar Z$ amounts then to add new symmetries, namely with respect to the medians of the triangles, corresponding to the root system $G_2$.  Let us describe this new polynomial system. Setting $s= Z+ \bar Z$ and $p= Z\bar Z$, formulae~\eqref{eq.deltoid} give
\beq\label{eq.G2}\bcas \Gamma(s,s)= p-s^2+s+1,\\
\Gamma(s,p)= s^2-2p-\frac{3}{2} sp + \frac{s}{2},\\
\Gamma(p,p)= s^3-3p^2-3sp+p\\
 \LL^{(\lambda )}(s)= -\lambda s, ~ \LL^{(\lambda)}(p)=1-(2\lambda+1)p
\ecas \eeq

Let us call $\tilde  \LL^{(\lambda)}$ this operator acting on functions (indeed polynomials) in the variables $(s,p)$.

From equation~\eqref{eq.G2}, it is clear that  the operator $\tilde \LL^{(\lambda)}$ preserves the set of polynomials in the variables $(s,p)$ (this just translates the invariance of $  \LL^{(\lambda)}$ under $Z\mapsto \bar Z$). However, because of the term $p^3$ in the coefficient $\Gamma(p,p)$, it does not preserve the degree of the polynomial. But things work better if we decide that the degree of $s^rp^t$ is $s+2t$. Then, with this new  notion of degree, a polynomial $Q(s,p)$ of degree $k$ is transformed under $\tilde \LL^{(\lambda)}$  into a polynomial of degree $k$. One may therefore find an orthonormal basis for $\tilde\  \LL^{(\lambda)}$ as polynomials in the variables $(s,p)$. In fact, it is nothing else than the symmetric polynomials $P_{n,k}(Z,\bar Z)$ expressed as polynomials in $(s,p)$.  Since now, the eigenspaces are one dimensional, one may play the same operation with the operator  $\hat L^{(\lambda)}$, but now with the projection $\Psi : Z\mapsto \frac{1}{2}(Z+\bar Z)= s, Z\bar Z= p$. We shall obtain the true hypergroup property for this polynomial family, through the projection $\Omega_1\mapsto R_\theta(z)+ \overline{R_\theta(z)}= R_\theta(z)+ R_{-\theta}(\bar z)$. 

 The image operator may be diagonalized in a family of orthogonal polynomials (in the variables $(s,p)$) for the image measure. This model    does not appear in the~\cite{BOZ2013} classification, since in this case the orthogonal polynomials must be ranked according to a degree which is $2\deg(p)+\deg(s)$, whereas in~\cite{BOZ2013}, the polynomials are ranked according to their usual degree.  The boundary of $\Omega_2$ is indeed the set where the determinant of the metric vanishes. This determinant may be written as 
$$\frac{1}{4}(s^2-4p)(3s^2+12sp+6p-4s^3-1),$$ and the boundary of $\Omega_2$ is a degree 5 algebraic curve. 
The curves $s^2-4p=0$ and $ 3p^2+12sp+6p-4s^3-1=0$ correspond respectively  to the images under $\Psi$ of the line $Z= \bar Z$ (the real axis), and of the boundary of $\Omega$ (the deltoid curve). It turns out that this last curve is a cuspidal cubic. Setting
$$s= x-1, p= y-2x+1,$$  (which corresponds to some affine change of coordinates), they are transformed in
$$4x^3-3y^2=0, ~x^2+6x-4y-3=0.$$ These two curves (cubic an parabola) cross  in the point $(1/3, -2/9)$, and are tangent to the second order at the point $(3,6)$ (in  the variables $(x,y)$).  The cuspidal point in the cubic $(0,0)$ is the image of $j$ and $\bar j$ in the deltoid (two of the cuspidal points of the deltoid), while the point $(3,6)$ is the image of the third cusp of the deltoid, that is the point $1$, which is also on the line $Z= \bar Z$.

Moreover, with $Q_1= s^2-4p$ and $Q_2= 3p^2+12sp+6p-4s^3-1$,  both $Q_1$ and $Q_2$ satisfy an equation similar to equation~\eqref{cond.poly.delt2}, and more precisely
$$\bcas \Gamma(\log Q_1,s)= -2s-2, ~
\Gamma(\log Q_1,p)= -3p-2s+1,\\
\Gamma(\log Q_2,s)= -3s,~
\Gamma(\log Q_2,p)= -6p
\ecas
$$
This shows that the operator defined through equations~\eqref{eq.G2} has reversible measure $C_\lambda Q_1^{-1/2}Q_2^{(2\lambda-5)/6} dsdp$, which  is of course the image of the measure $\mu^{(\lambda)}$ through the projection  $\Psi$. This operator therefore satisfies the usual hypergroup property, with reference point the image of $1$, which is $(2,1)$, that is the point where the cuspidal cubic and the parabola are bi-tangent to each other.

But from the general presentation of~\cite{BOZ2013}, there is now a two-parameter family of measures, namely $\mu_{\alpha, \beta}(ds,dp)= C_{\alpha, \beta} Q_1^{\alpha_1}Q_2^{\alpha_2} dsdp$ on this set $\Omega_2$, for which there exist a family of orthogonal polynomials which are eigenvectors of a diffusion operator.

 The conditions under which those measures are finite are   $\alpha_1>-1$, $\alpha_2>-5/6$,  and $\alpha_1+\alpha_2> -4/3$. The condition $\alpha_2>-5/6$ and $\alpha_1+\alpha_2> -4/3$ correspond respectively to the integrability conditions around the cusp of the cubic and the double tangent point. 
For the double tangent point,  to check the integrability condition, one may reduce, up to an affine transformation of the plane,  to check the integration condition for $$\int_0^1\int_0^2 \un_{y^2<x< y^2+ cy^3} (x-y^2)^\alpha_1(y^2+cy^3-x)^\alpha_2 dxdy.$$
 After a few change of variables,  this reduces to check the integrability condition for $\int_0^1 y^{3(\alpha_1+\alpha_2+1)} dy$.  The cusp is easier to deal with and quite immediate.

For these measures $\mu_{\alpha_1, \alpha_2}$, we  have  an associated operator $\LL_{\alpha_1, \alpha_2}$ (sharing the same $\Gamma$ operator, and defined through formula~\eqref{ipp}) and associated orthogonal polynomials  which are  eigenvectors of $\LL_{\alpha_1, \alpha_2}$.

We already  mentioned that this model reflects in fact the symmetries of the root system $G_2$.  In this context, the Weyl group acting on the roots has two orbits, corresponding to the two irreducible factors of the boundary of $\Omega_1$, and to the two parameters in the choice of the measure.  When $\alpha_1=-1/2$, the model is the direct image of the deltoid model through the projection $\Psi$, and the eigenvectors of the associated operators are just the rewriting of $P_{n,k}$ as polynomials in the variables $s=Z+ \bar Z$ and $p= Z\bar Z$.  The hypergroup property for this model is the direct consequence of the previous Section~\ref{sec.hgp.deltoid1}, but may also be reproved directly using the 6 dimensional model of Section~\ref{6.dim.model}.

 The case where $\alpha_1\neq -1/2$  remains open. Let us show however how to deal with a different value of $\alpha_1$, whenever $\alpha_2=-1/2$. This relies on a very specific property of this $G_2$ polynomial model.  We first observe that the domain $\Omega_2$ is stable under the transformation 
 $$\Psi_1 :(s,p)\mapsto (S= 3p-1,P=  1+s^3-3ps-6p).$$ Under this transformation, the parabolic part of the boundary is mapped onto the cubic one, and conversely. The transformation $\Psi_1$ comes in fact from the invariance $Z\mapsto jZ$ in the deltoid model, where instead of looking at functions of $(Z\bar Z, Z+\bar Z)$, we looked at functions of $(Z\bar Z, Z^3+ \bar Z^3)$, under a slight change of coordinates such that the domain $\Omega_2$ is invariant under the transformation.
 
 Moreover, one may check that, when $\alpha_1=-1/2$, the image of the operator $\LL_{-1/2, \alpha_2}$ is $\frac{1}{3}\LL_{\alpha_2, -1/2}$. It may be easily checked looking at $\Gamma(S,S), \Gamma(S,P), \Gamma(P,P), \LL_{-1/2, \alpha_2}(S)$ and $\LL_{-1/2, \alpha_2}(P)$. There is nothing similar when $\alpha_1\neq -1/2$. The transformation $\Psi_1$ is not a diffeomorphism, and the operators $\LL_{\alpha_1, \alpha_2}$ in general do not have images under $\Psi_1$ when $\alpha_1\neq -1/2$. From the previous scheme, we may now conclude to the hypergroup property for the orthogonal polynomials associated with $\LL_{\alpha, -1/2}$ when $\alpha\geq 0$.

 \bibliographystyle{amsplain}   
\bibliography{bib.bakry.Hypergr.deltoid}
\end{document}